\RequirePackage{fix-cm}
\documentclass[smallextended,envcountsame,envcountsect]{svjour3}
\smartqed
\usepackage{amsmath,amssymb,booktabs,tikz}
\IfFileExists{pgflibraryarrows.meta.code.tex}{%
  \usetikzlibrary{arrows.meta}%
  \tikzset{cycle arrows/.style={>=Stealth}}%
}{%
  \tikzset{cycle arrows/.style={>=stealth}}%
}
\usepackage{hyperref}
\hypersetup{colorlinks=false,pdfborder={0 0 0}}
\hypersetup{pdftitle={Minimum-order realizations of prescribed cycle lattices in two-letter digraphs},pdfauthor={Xuankai Rong, Zhilan Miao, Xinyao Li, Panwang Wang}}
\journalname{Noname}
\newcommand{\Z}{\mathbb Z}
\newcommand{\N}{\mathbb N}
\newcommand{\R}{\mathbb R}
\newcommand{\E}{\mathbb E}
\newcommand{\Var}{\operatorname{Var}}
\newcommand{\rank}{\operatorname{rank}}
\newcommand{\mc}{M_{\mathrm{comp}}}
\newcommand{\mn}{M_{\mathrm{nd}}}

\newcommand{\dd}{\,\mathrm d}
\newcommand{\rdoi}[1]{\url{https://doi.org/#1}}
\begin{document}
\title{Minimum-order realizations of prescribed cycle lattices in two-letter digraphs}
\titlerunning{Minimum-order realizations of cycle lattices}
\author{Xuankai~Rong \and Zhilan~Miao \and Xinyao~Li \and Panwang~Wang\textsuperscript{*}}
\authorrunning{X. Rong et al.}
\institute{Xuankai Rong and Zhilan Miao \at
School of Information Science and Technology, Shijiazhuang Tiedao University, Shijiazhuang, China\\
\email{rxuan894@gmail.com} (Xuankai Rong)\\
\email{499325989@qq.com} (Zhilan Miao)
\and Xinyao Li \at
School of Computer and Communication Engineering, University of Science and Technology Beijing, Beijing, China\\
\email{starrysea0531@outlook.com}
\and Panwang Wang\textsuperscript{*} (Corresponding author) \at
Department of Mathematics and Physics, Shijiazhuang Tiedao University, Shijiazhuang, China\\
\email{wangpw@stdu.edu.cn}}
\date{}
\maketitle
\begin{abstract}
Let \(L\le\Z^2\) be a sublattice of finite index. Among strongly connected, partial deterministic edge-labelled digraphs over a binary alphabet, let \(M(L)\) denote the minimum number of vertices required for the letter-count vectors (Parikh vectors) of closed directed walks to generate \(L\). We prove that \(M(L)\) equals one less than the minimum semiperimeter of an axis-aligned rectangle anchored at the origin that contains two linearly independent points of \(L\). This minimum is attained by a partial deterministic directed theta graph and remains sharp even over unrestricted (nondeterministic) realizations. Moreover, any realization containing \(k\ge2\) pairwise vertex-disjoint directed cycles of length at least \(g\ge1\) whose count vectors span \(\R^2\) requires at least \(1+(k-2)g\) vertices beyond \(M(L)\); this excess bound is sharp for all \(k\ge2\) and \(g\ge1\). We further determine the exact minimum overlap slack between two basis cycles with a prescribed number of common components, and establish an additive formula for direct sums over disjoint alphabet blocks. Finally, for a uniformly random nonzero congruence slope of prime index \(p\), the normalized minimum \(M(L)/\sqrt p\) converges in distribution to the minimum rectangle semiperimeter of a Haar-distributed unimodular lattice. Its expectation converges to the finite mean of that functional, and the unnormalized variance satisfies \(\Var(M(L))=(12/\pi^2)p\log p+O(p)\).

\end{abstract}

\section{Introduction}

\paragraph{The realization problem.}
The Parikh map records the number of occurrences of each letter in a word \cite{parikh}. Reading the edge labels on a directed walk gives such a word. This paper addresses an inverse structural question: given the lattice generated by closed-walk counts, how few vertices suffice to realize it?

Let \(G\) be a finite, strongly connected directed multigraph with at least one edge, each edge carrying a label from a finite alphabet. Write \(V(G)\) for its vertex set and \(|V(G)|\) for its order. Fix the alphabet \(\Sigma=\{a,b\}\) in that order. Let \(\Sigma^*\) denote all finite words, including the empty word, and let \(\N=\{0,1,\ldots\}\). The Parikh map \(\Psi:\Sigma^*\to\N^2\) sends \(w\) to \((|w|_a,|w|_b)\), where \(|w|_i\) counts occurrences of letter \(i\). For a walk \(W\), write \(\Psi(W)\) for the vector of its label word, counting each edge traversal. Thus a path reading \(aab\) has vector \((2,1)\), and \(\Psi(PQ)=\Psi(P)+\Psi(Q)\) for consecutive walks \(P,Q\). A letter exponent denotes repetition. For larger alphabets use one coordinate per letter. Define the \emph{cycle lattice} of \(G\) by
\[
 L(G)=\langle\Psi(W):W\text{ is a closed directed walk in }G\rangle_{\Z},
\]
where the angle brackets denote the \(\Z\)-linear span. A point of \(L(G)\) need not itself be realized as the count of a single closed walk. Equivalently, \(L(G)\) is the image of the integral flow lattice of \(G\) under the map summing edge weights by label. Indeed, strong connectivity yields a strictly positive integral circulation (obtained by summing closed walks covering all edges); adding a sufficiently large multiple of this circulation to any integral flow produces a nonnegative circulation, which decomposes into directed cycles, and subtracting back gives the claim.

We say that \(G\) \emph{realizes} \(L\) if \(L(G)=L\). A labelled digraph is \emph{partial deterministic} if every vertex has at most one outgoing edge of each label, and \emph{complete deterministic} if every vertex has exactly one. Let \(M(L)\), \(\mc(L)\), and \(\mn(L)\) denote the minimum number of vertices of a realization of \(L\) in the partial deterministic, complete deterministic, and unrestricted classes, respectively, with the convention that the minimum over an empty class is \(+\infty\). Unless alphabet blocks are specified, all graphs considered below are binary.

\paragraph{Related work.}
Our setting sits at the intersection of automata theory, structural graph theory, and the geometry of numbers. In formal language theory, Champarnaud et al.\ \cite{geometry} use simple-cycle counts in strongly connected automata to construct linear Diophantine systems in the study of geometrical regular languages. A broader graph-theoretic analogue is the lattice of integral flows studied by Bacher et al.\ \cite{flows}, where an integral flow assigns an integer to each oriented edge subject to flow conservation at every vertex; summing edge weights by label projects this high-dimensional lattice onto our cycle lattice.

Two established lines of work examine how letter counts govern automaton size and expressiveness. Parikh automata, introduced by Klaedtke and Rue\ss{} \cite{klaedtke} and developed algebraically by Cadilhac et al.\ \cite{cadilhac}, accept a word whenever a successful run accumulates a total weight in a prescribed semilinear set (a finite union of translates of finitely generated additive monoids). Lavado et al.\ \cite{lavado} bound the state complexity of deterministic automata preserving the entire Parikh image of a language recognized nondeterministically, while Hoffmann \cite{hoffmann} investigates commutative closure for languages accepted by permutation automata. Unlike these models, our digraphs carry no designated initial or accepting states: we prescribe only the integer subgroup generated by closed-walk counts, so distinct realizations of the same lattice may exhibit different semigroups of attainable return counts.

From a combinatorial and geometric perspective, Echols et al.\ \cite{orbits} characterize the realizable pairs of lengths and step counts---and enumerate primitive periodic orbits---in \(2\)-regular circulant digraphs, where the lattice of relations between the two generators plays a central role. Cai et al.\ \cite{routing} apply congruence lattices to estimate shortest loop lengths in circulant graphs, and Marklof and Str\"ombergsson \cite[Theorem~1]{diameters} establish a limit law for diameters of random directed circulant graphs via the covering radius of a simplex for a random unimodular lattice. On a Riemannian torus, Martelli et al.\ \cite[Theorem~4.2]{spines} prove that a spine of minimal length has theta topology, where edges are geodesic arcs and the objective is Riemannian length rather than discrete vertex count.

Finally, for two oriented simple closed curves on a torus, the algebraic intersection number equals the determinant of their homology coordinates in a positively oriented basis \cite[p.~29]{farb}. For a graph embedded on a torus, Schrijver \cite[Theorem~1(i), p.~147]{schrijver} lower-bounds the maximum number of disjoint noncontractible circuits in terms of the face width, defined as the minimum of \(\lvert C\cap G\rvert\) over all noncontractible closed curves \(C\) on the surface. In our setting, a full-rank lattice \(L\le\Z^2\) defines the torus \(\R^2/L\) with grid phases in the quotient group \(\Z^2/L\), whose elements are the cosets \(z+L=\{z+v:v\in L\}\) for \(z\in\Z^2\); two cycles that are vertex-disjoint in the original digraph may nevertheless collide at a common phase.

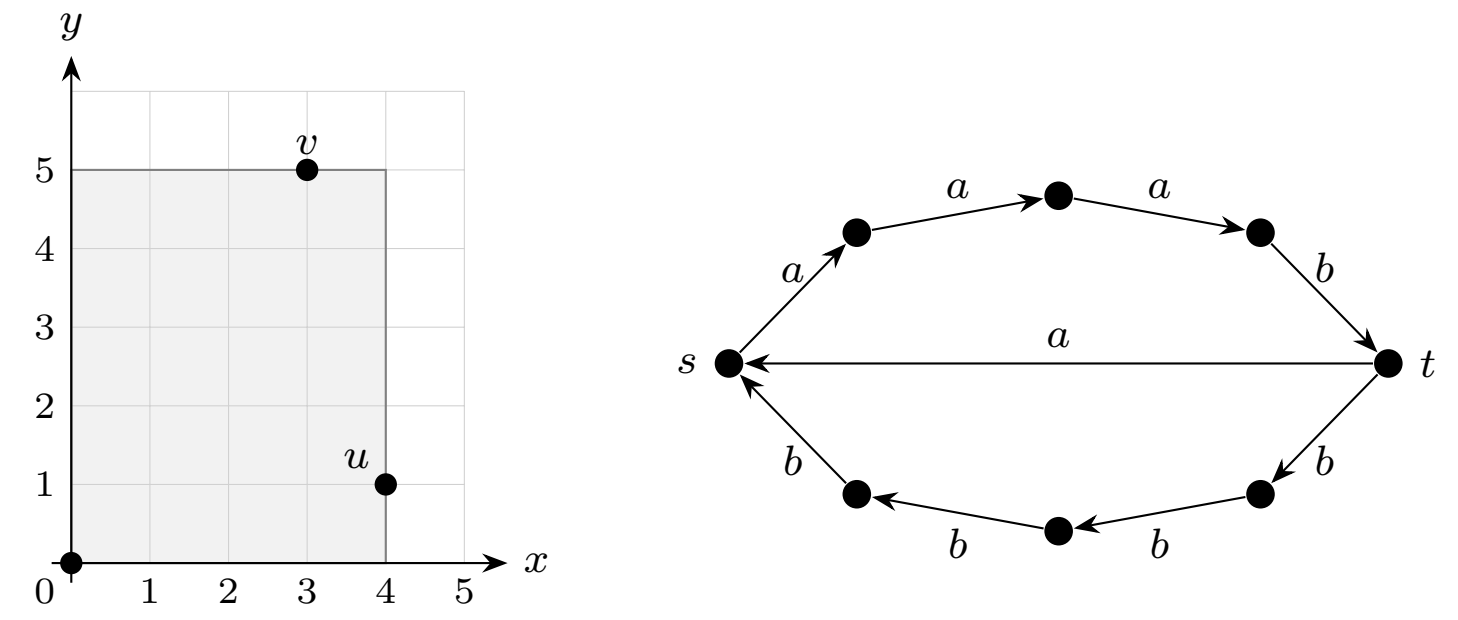
\begin{figure}[tb]
\centering
\begin{tikzpicture}[cycle arrows,font=\small]
\begin{scope}[x=.52cm,y=.52cm]
\fill[black!5] (0,0) rectangle (4,5);
\draw[step=1,black!18,very thin] (0,0) grid (5,6);
\draw[black!50] (0,0) rectangle (4,5);
\draw[->] (-.25,0)--(5.55,0) node[right] {$x$};
\draw[->] (0,-.25)--(0,6.45) node[above] {$y$};
\foreach \x in {1,2,3,4,5}{\node[below,font=\scriptsize] at (\x,0) {\x};}
\foreach \y in {1,2,3,4,5}{\node[left,font=\scriptsize] at (0,\y) {\y};}
\foreach \x/\y in {0/0,4/1,3/5}{\fill (\x,\y) circle (2.1pt);}
\node[above left] at (4,1) {$u$}; \node[above] at (3,5) {$v$};
\node[below left,font=\scriptsize] at (0,0) {$0$};
\end{scope}
\begin{scope}[xshift=4.35cm,yshift=1.32cm,x=.89cm,y=.94cm,
  dot/.style={circle,fill,inner sep=1.9pt}]
\node[dot,label=left:$s$] (s) at (0,0) {};
\node[dot,label=right:$t$] (t) at (4.9,0) {};
\foreach \id/\x/\y in {u/.95/.92,v/2.45/1.18,w/3.95/.92,i/3.95/-.92,j/2.45/-1.18,k/.95/-.92}{\node[dot] (\id) at (\x,\y) {};}
\foreach \a/\b/\sym in {s/u/a,u/v/a,v/w/a,w/t/b}{\draw[->] (\a)--node[above]{$\sym$}(\b);}
\draw[->] (t)--node[above]{$a$}(s);
\foreach \a/\b in {t/i,i/j,j/k,k/s}{\draw[->] (\a)--node[below]{$b$}(\b);}
\end{scope}
\end{tikzpicture}
\caption{Two descriptions of the optimum for \(x\equiv4y\pmod {17}\). The shaded rectangle contains \(u=(4,1)\) and \(v=(3,5)\). In the graph, the upper path from \(s\) to \(t\) reads \(a^3b\); returning via the middle edge gives \(u\), and returning along the lower path gives \(v\). The semiperimeter is nine and the graph has eight vertices.}\label{fig:theta}
\end{figure}

As a concrete illustration, tracking \(x-4y\bmod 17\) yields a \(17\)-vertex complete deterministic realization of the lattice \(x\equiv 4y\pmod{17}\). Yet the two return words \(a^3ba\) and \(a^3b^5\) have counts generating the same lattice and can share their common prefix \(a^3b\), producing the \(8\)-vertex partial deterministic graph in Fig.~\ref{fig:theta}. A rectangle of semiperimeter at most eight has area at most sixteen, so two points in it with nonnegative coordinates cannot have determinant of absolute value at least seventeen. Theorem~\ref{th:main} therefore proves optimality. By contrast, changing the congruence to \(x\equiv y\pmod{17}\) forces all \(17\) vertices by Proposition~\ref{pr:extreme}.

\paragraph{Overview of results.}
In Section~\ref{sec:rectangle}, we establish the exact rectangle formula \(\mn(L)=M(L)=B(L)-1\), where \(B(L)\) is the least semiperimeter of a first-quadrant rectangle anchored at the origin containing two independent lattice points (Theorem~\ref{th:main}), showing that nondeterminism offers no vertex savings and that an optimal partial deterministic theta graph always exists. We also present a self-contained alternative proof of the lower bound working directly in the original graph, as well as the rank-one case. Section~\ref{sec:bounds} derives sharp index bounds (Corollary~\ref{co:index}), shows that \(\mc(L)=D\), where \(D=[\Z^2:L]\) (Proposition~\ref{pr:complete}), proves that \(M(L)\) is computable in time polynomial in the binary encoding length of a Hermite basis (Corollary~\ref{co:algorithm}), and establishes an exact additive formula for direct sums over disjoint alphabet blocks (Theorem~\ref{th:blocks}).

Section~\ref{sec:packing} investigates structural constraints on cycles inside realizations. In every minimum realization, two cycles with linearly independent counts must share a vertex; Theorem~\ref{th:packing} determines the exact minimum vertex excess \(1+(k-2)g\), over all lattices and realizations, forced by \(k\) pairwise vertex-disjoint cycles of length at least \(g\) whose counts span \(\R^2\). Theorem~\ref{th:slack} then establishes the exact minimum overlap slack between two basis cycles whose intersection has \(c\ge2\) connected components. Finally, Section~\ref{sec:random} analyzes random congruence lattices \(x\equiv sy\pmod p\) of prime index \(p\), characterizing the extremal slopes (Proposition~\ref{pr:extreme}) and proving both the leading variance asymptotics (Theorem~\ref{th:random}) and the distributional and mean limits under the Haar measure on \(\mathrm{SL}_2(\R)/\mathrm{SL}_2(\Z)\) (Theorem~\ref{th:limit}). Section~\ref{sec:concluding} closes with open problems.

\section{An exact rectangle formula}\label{sec:rectangle}
Write \(u=(u_x,u_y)\) for a vector in \(\Z^2\). For \(u,v\in\Z^2\), set \(\det(u,v)=u_xv_y-u_yv_x\); we call \(u,v\) \emph{independent} if \(\det(u,v)\ne0\). A full-rank sublattice \(L\le\Z^2\) contains an independent pair, and its index \(D=[\Z^2:L]\) is the number of cosets of \(L\) in \(\Z^2\). We write \(\|u\|_1=|u_x|+|u_y|\) and \(\|u\|_\infty=\max(|u_x|,|u_y|)\). A basis \(u,v\in\N^2\) of \(L\) is \emph{oppositely ordered} if, after possibly swapping its elements, \(u_x>v_x\) and \(u_y<v_y\). For a full-rank sublattice \(L\le\Z^2\), define
\begin{align}
 B(L)&=\min\{w+h:\ w,h\in\N,\
   \rank\langle L\cap([0,w]\times[0,h])\rangle_{\Z}=2\},\label{eq:box}\\
 F(L)&=\min\{u_x+v_y-1:\ (u,v)\in\mathcal A(L)\}.\nonumber
\end{align}
Here \(\mathcal A(L)\) consists of all ordered bases \((u,v)\) of \(L\) with \(u,v\in\N^2\), \(u_x>v_x\) and \(u_y<v_y\). In \eqref{eq:box}, the rank-two condition simply states that the closed rectangle \([0,w]\times[0,h]\) contains two independent points of \(L\); we call \(w+h\) its \emph{semiperimeter}. Notice that a Hermite basis \((h_{11},h_{21}),(0,h_{22})\) with \(h_{11},h_{22}>0\) and \(0\le h_{21}<h_{22}\) is admissible (and oppositely ordered), with \(D=h_{11}h_{22}\).

In the plane, Pick's theorem implies that a lattice triangle containing no lattice points other than its three vertices has area \(1/2\) in lattice coordinates \cite[Theorem~2.8, p.~41]{beck}. The elementary determinant argument below establishes the basis step directly while keeping the oppositely ordered pair inside the prescribed rectangle.

\begin{lemma}\label{le:box}
Every closed rectangle \(R=[0,A]\times[0,B]\) containing two independent points of \(L\) contains an oppositely ordered basis of \(L\). Consequently, \(F(L)=B(L)-1\).
\end{lemma}
\begin{proof}
Among all independent pairs \(u,v\in L\cap R\), choose one minimizing \(|\det(u,v)|\). Suppose for a contradiction that \(u,v\) do not generate \(L\). Then their half-open fundamental parallelogram contains a representative \(z=\alpha u+\beta v\) of a nonzero coset, where \(0\le\alpha,\beta<1\). If \(\alpha+\beta>1\), replace \(z\) by \(u+v-z\) and its coefficients by \(1-\alpha,1-\beta\). We then have \(0<\alpha+\beta\le1\) with both coefficients still strictly less than \(1\), so \(z\) lies in the triangle with vertices \(0,u,v\), and hence in \(R\). If \(\alpha>0\), then \(\det(z,v)=\alpha\det(u,v)\); otherwise \(\beta>0\) and \(\det(u,z)=\beta\det(u,v)\). In either case, replacing one vector by \(z\) produces an independent pair in \(L\cap R\) with strictly smaller absolute determinant, a contradiction. Thus \(u,v\) form a basis of \(L\).

If one of \(u,v\) dominates the other coordinatewise, subtract the smaller vector from the larger. Both vectors remain in \(R\), the generated lattice is unchanged, and the sum of all four coordinates strictly decreases. Repeating this step must therefore terminate with an oppositely ordered basis of \(L\) contained in \(R\). Applying this to a minimizing rectangle yields \(F(L)\le B(L)-1\), while the bounding rectangle of a minimizing oppositely ordered basis gives the reverse inequality. \qed
\end{proof}

\begin{theorem}\label{th:main}
For every full-rank sublattice \(L\le\Z^2\),
\[
 \mn(L)=M(L)=F(L)=B(L)-1.
\]
Moreover, this common minimum is attained by a partial deterministic directed theta graph, consisting of an initial path joining two branch vertices together with two return paths whose interiors are pairwise disjoint from each other and from the initial path (with the two branch vertices identified when the initial path is empty).
\end{theorem}

\begin{figure}[tb]
\centering
\begin{tikzpicture}[cycle arrows,font=\small,
 first/.style={black,thick,->},second/.style={black!75,thick,dashed,->},
 state/.style={draw,circle,fill=white,minimum size=5.3mm,inner sep=0pt}]
\begin{scope}[x=.85cm,y=.72cm]
\draw[step=1,black!18,very thin] (0,0) grid (3,3);
\draw[->,black!55] (-.15,0)--(3.4,0) node[right] {$x$};
\draw[->,black!55] (0,-.18)--(0,3.6) node[above] {$y$};
\draw[first] (0,1)--(1,1); \draw[first] (1,1)--(2,1); \draw[first] (2,1)--(2,2);
\draw[second] (1,0)--(1,1); \draw[second] (1,1)--(1,2);
\draw[second] (1,2)--(1,3); \draw[second] (1,3)--(2,3);
\fill (1,1) circle (2.1pt);
\node[below left,font=\scriptsize] at (1,1) {$(1,1)$};
\node[right] at (2,2) {$C$}; \node[right] at (2,3) {$C'$};
\end{scope}
\begin{scope}[xshift=5.35cm,yshift=1.45cm]
\node[state] (f) at (0,0) {$4$};
\node[state] (z) at (3.05,0) {$0$};
\node[state] (h) at (1.52,1.02) {$3$};
\node[state] (e) at (1.52,-1.02) {$2$};
\draw[first] (f)--node[above]{$a$}(z);
\draw[first] (z)--node[above right]{$b$}(h);
\draw[first] (h)--node[above left]{$a$}(f);
\draw[second] (f)--node[below left]{$b$}(e);
\draw[second] (e)--node[below right]{$b$}(z);
\end{scope}
\end{tikzpicture}
\caption{Grid crossing and ear extraction in the quotient with phase \(x-2y\bmod5\). Left: one lifted lap of \(C\) (solid, \(aab\)) meets one of \(C'\) (dashed, \(bbba\)) at \((1,1)\), of phase \(4\). Right: retain \(C\) and the ear \(4\to2\to0\) of \(C'\). Closing the ear along \(0\to3\to4\) yields count \((1,3)\), independent of the count \((2,1)\) of \(C\). The retained graph has five edges and four vertices.}\label{fig:ear}
\end{figure}
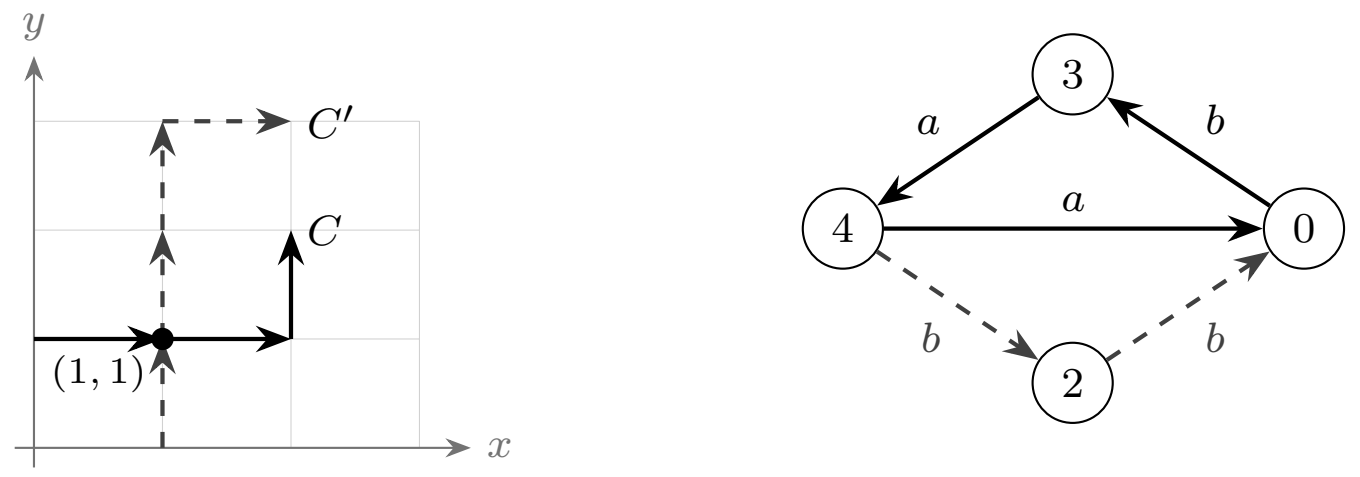

\begin{proof}
\emph{Construction.} Fix an oppositely ordered basis \(u,v\) attaining \(F(L)\). Construct a directed path from \(s\) to \(t\) labelled by \(a^{v_x}b^{u_y}\), together with two return paths from \(t\) to \(s\) labelled by \(a^{u_x-v_x}\) and \(b^{v_y-u_y}\), respectively, taking all internal vertices distinct. When \(v_x=u_y=0\), the initial path is empty and we identify \(s=t\). Because \(u_x-v_x>0\) and \(v_y-u_y>0\), the two outgoing edges at \(t\) carry distinct labels, so the resulting digraph is partial deterministic. Every vertex lies on a directed cycle through \(s\), and the graph has \(u_x+v_y\) edges and \(u_x+v_y-1\) vertices. Every nonempty closed directed walk visits \(s\); after a cyclic shift, decomposing at successive visits to \(s\) expresses it as a concatenation of the two basis cycles. Hence the cycle lattice is precisely \(L\), giving \(\mn(L)\le M(L)\le F(L)\).

\emph{Passage to the quotient.} Let \(G\) be any realization of \(L\), possibly nondeterministic. Fix a root vertex \(o\in V(G)\) and assign to each vertex \(q\in V(G)\) the \emph{phase} \(\phi(q)=\Psi(P)+L\in\Z^2/L\). Here \(P\) is any directed path from \(o\) to \(q\): \(\Psi(P)=(|P|_a,|P|_b)\) counts its edge labels, and \(\phi(q)\) is the coset of that vector modulo \(L\). Appending a common return path from \(q\) to \(o\) to any two choices of \(P\) shows that their count vectors differ by an element of \(L\), so \(\phi\) is well defined. Form the quotient digraph \(H\) whose vertex set is \(\phi(V(G))\) and whose edges are the images of the labelled edges of \(G\), retaining at most one copy of each labelled edge between any ordered pair of phases. Let \(e_a=(1,0)\) and \(e_b=(0,1)\). From any given phase, every edge labelled \(i\in\{a,b\}\) translates the phase by \(e_i\), so \(H\) is partial deterministic and strongly connected. Every closed walk in \(G\) projects to a closed walk in \(H\) with the same count, while every closed walk in \(H\) has zero net displacement in \(\Z^2/L\) and hence count in \(L\). Thus
\[
 L(H)=L,\qquad |V(H)|\le |V(G)|.
\]

\emph{Intersection of cycles.} Since closed directed walks decompose into simple directed cycles, \(H\) contains two simple cycles \(C,C'\) with independent counts \(u,v\). Choose integer lifts of their base vertices and lift each edge of label \(a\) or \(b\) to a unit east or north step in \(\Z^2\). Extending periodically in both directions yields two doubly infinite lattice paths, each having a unique vertex at every integer time \(t=x+y\). Let \(X(t)\) and \(Y(t)\) denote their horizontal coordinates at time \(t\), and put \(h=\|u\|_1\) and \(k=\|v\|_1\) for their periods in time. Over \(hk\) steps, the first path completes \(k\) laps and the second completes \(h\), so the integer-valued function \(f=X-Y\) satisfies
\[
 f(t+hk)-f(t)=ku_x-hv_x=\det(u,v)\ne0.
\]
Consequently, \(f\) attains both positive and negative values along multiples of \(hk\). Since \(f\) changes by at most \(1\) at each unit time step, \(f(t_0)=0\) for some integer \(t_0\). At time \(t_0\), the two lifts share both horizontal and vertical coordinates, so \(C\) and \(C'\) intersect at a common vertex of \(H\), as illustrated in Fig.~\ref{fig:ear}.

\emph{Extraction of an ear.} Here an \emph{ear} is a positive-length directed walk with endpoints on \(C\) and distinct internal vertices outside \(C\); its endpoints may coincide. Decompose \(C'\) at its visits to \(C\) in cyclic order along \(C'\), yielding segments \(Q_i\colon x_i\to x_{i+1}\) (with cyclic indices), and for each \(i\) let \(R_i\) be the forward segment along \(C\) from \(x_{i+1}\) to \(x_i\) (taking \(R_i\) empty when \(C\) and \(C'\) meet at a single vertex). Because each common vertex appears equally often as the start and the end of a return segment, the formal sum of the \(R_i\) forms a nonnegative circulation on the simple cycle \(C\), so every edge of \(C\) appears with the same multiplicity \(m\in\N\). Thus
\[
 \sum_i\Psi(Q_iR_i)=v+mu.
\]
Since \(\det(u,v+mu)=\det(u,v)\ne0\), at least one summand \(\Psi(Q_iR_i)\) is independent of \(u\). Fix such an index \(i\). Because the interior of \(Q_i\) is disjoint from \(C\), if \(Q_i\) shared an edge with \(C\) it would consist solely of that edge, so \(Q_iR_i\) would be a full lap of \(C\), giving \(\Psi(Q_iR_i)=u\), contradicting the choice of \(i\). Therefore \(K=C\cup Q_i\) has exactly one more edge than vertex: an ear of length \(\ell\) contributes \(\ell\) new edges and \(\ell-1\) new vertices (a count that remains valid when \(Q_i\) is a loop or an edge parallel to an edge of \(C\)).

Let \(A\) and \(B\) be the numbers of edges of \(K\) labelled \(a\) and \(b\), respectively. Both \(C\) and \(Q_iR_i\) traverse each edge of \(K\) at most once, so their independent count vectors lie in \(L\cap([0,A]\times[0,B])\). Hence
\[
 B(L)\le A+B=|V(K)|+1\le |V(H)|+1\le |V(G)|+1.
\]
Combining this lower bound with the upper bound from the construction and Lemma~\ref{le:box} completes the proof. \qed
\end{proof}

\paragraph{Alternative proof of the lower bound.}
The Lean formalization of Theorem~\ref{th:main} establishes the lower bound directly in the original graph \(G\), without passing to the phase quotient. Fix a simple cycle \(C\) in \(G\) with count \(u\). Assign to each vertex \(x\in V(C)\) the height \(h(x)=\det(u,z)\), where \(z\) is the count of the prefix of \(C\) from a fixed root to \(x\); a full traversal of \(C\) adds \(\det(u,u)=0\), so \(h\) is well defined on \(V(C)\). Extend \(h\) arbitrarily to \(V(G)\setminus V(C)\), and assign to every directed path from \(x\) to \(y\) with count \(v\) the adjusted weight \(\det(u,v)+h(x)-h(y)\). These weights are additive under concatenation, and every forward segment of \(C\) has adjusted weight zero. Because \(L(G)\) has full rank, there exists a closed walk with count independent of \(u\); choose paths from the root to its base and back. The two root returns obtained with and without the inserted walk have different weights, so at least one has nonzero weight.

Choose a shortest walk of nonzero weight among all walks with both endpoints on \(C\) and all closed walks disjoint from \(C\). A walk of the first kind cannot visit \(C\) internally, for splitting at an internal visit would yield a strictly shorter candidate of nonzero weight; similarly, no exterior vertex can repeat, as detaching the resulting closed subwalk would leave either the subwalk or the remaining walk with nonzero weight. Thus the minimizer is either a simple cycle disjoint from \(C\) or an ear with distinct internal vertices outside \(C\). In the ear case, closing the ear along \(C\) produces a closed walk whose adjusted weight equals the determinant of \(u\) with its count vector. The coordinatewise sum of the counts of \(C\) and the ear therefore defines a rectangle containing two independent return counts, of semiperimeter at most \(|V(G)|+1\) (since an ear of length \(\ell\) uses \(\ell-1\) exterior vertices). In the disjoint case, the sum of the two cycle lengths is already at most \(|V(G)|\).

\paragraph{Rank-one lattices.}
For a rank-one sublattice \(L\le\Z^2\), the real quotient \(\R^2/L\) is a cylinder; note that strong connectivity together with the presence of at least one edge precludes rank zero. If \(L=\Z w\), then \(L\) admits a realization if and only if \(w\) or \(-w\) lies in \(\N^2\). When this holds, a single directed cycle realizes \(L\) with \(\|w\|_1\) vertices. Choosing the sign of \(w\) so that \(w\in\N^2\), every nonempty simple cycle in any realization has count equal to a positive integer multiple of \(w\), establishing minimality even over unrestricted digraphs. Note that the generator \(w\) need not be primitive in \(\Z^2\).

\section{Index bounds, computation and alphabet blocks}\label{sec:bounds}
Write \(\lfloor x\rfloor\) for the greatest integer at most \(x\), and \(\lceil x\rceil\) for the least integer at least \(x\).

\begin{corollary}\label{co:index}
For \(L\) in the Hermite form of Section~\ref{sec:rectangle} and \(\nu=\gcd(h_{21},h_{22})\),
\[
 \max\{2\sqrt D-1,\,h_{11}+\nu-1\}\le M(L)\le h_{11}+h_{22}-1.
\]
In particular, \(M(\alpha\Z\times\beta\Z)=\alpha+\beta-1\) for all positive integers \(\alpha,\beta\). The largest finite index of a cycle lattice of a strongly connected \(n\)-vertex binary digraph is \(\lfloor(n+1)^2/4\rfloor\).
\end{corollary}
\begin{proof}
The Hermite basis is oppositely ordered and yields the upper bound. For an optimal oppositely ordered basis \(u,v\), we have \(D=u_xv_y-u_yv_x\le u_xv_y\), while \(u_x>0\) is divisible by \(h_{11}\) and \(v_y>0\) is divisible by \(\nu\). The lower bounds follow from the arithmetic--geometric mean inequality and divisibility. The divisibility lower bound equals the Hermite upper bound for rectangular lattices \(\alpha\Z\times\beta\Z\); choosing \(\alpha=\lfloor(n+1)/2\rfloor\) and \(\beta=\lceil(n+1)/2\rceil\) attains the maximum index for \(n\) vertices. \qed
\end{proof}

\begin{proposition}\label{pr:complete}
For every full-rank sublattice \(L\le\Z^2\), \(\mc(L)=D\).
\end{proposition}
\begin{proof}
In a complete deterministic realization, the phase image contains \(0\) and is closed under translation by both coordinate unit vectors \(e_a,e_b\). Since \(\Z^2/L\) is finite, these translations generate the group even as a monoid (as the inverse of each generator is a nonnegative multiple of itself), so the phase image has cardinality \(D\). Conversely, the Cayley digraph on \(\Z^2/L\) with generators \(e_a,e_b\) is complete deterministic and strongly connected. Its closed-walk counts are precisely \(L\cap\N^2\), which generate \(L\) because every \(w\in L\) can be written as \((w+k(D,D))-k(D,D)\) with both terms in \(L\cap\N^2\) for sufficiently large \(k\). \qed
\end{proof}

For the language of all words whose Parikh vector lies in \(L\), these \(D\) phases coincide with the Myhill--Nerode equivalence classes \cite{nerode}, as a suffix cancelling one phase distinguishes it from every other. Proposition~\ref{pr:complete} shows that the same state count persists under the closed-walk specification without accepting states.

For positive \(g\), we write \(f=O(g)\) if \(|f|\le Cg\) eventually for some constant \(C>0\), and \(f=\Theta(g)\) if \(cg\le f\le Cg\) eventually for some constants \(c,C>0\).

\begin{corollary}\label{co:algorithm}
Let \(H=\left(\begin{smallmatrix}h_{11}&0\\ h_{21}&h_{22}\end{smallmatrix}\right)\) be the Hermite basis matrix of Section~\ref{sec:rectangle}, so that \(L=H\Z^2\). With the entries of \(H\) encoded in binary, the value \(M(L)\) is computable in time polynomial in \(\log(2D)\).
\end{corollary}
\begin{proof}
Vector inequalities are coordinatewise. Given a budget \(m\), introduce two column vectors \(r,t\in\Z^2\) (four integer variables), set \(u=Hr\) and \(v=Ht\), and impose the linear inequalities
\[
 u,v\ge0,\qquad u_x-v_x\ge1,\qquad v_y-u_y\ge1,\qquad u_x+v_y\le m+1.
\]
Any feasible pair \((u,v)\) is linearly independent and lies in a rectangle of semiperimeter at most \(m+1\), whence \(M(L)\le m\) by Theorem~\ref{th:main}. Conversely, an optimal oppositely ordered basis is feasible whenever \(M(L)\le m\). This integer linear program has fixed dimension \(4\), and for \(1\le m\le D\) every coefficient has bit-length \(O(\log(2D))\). Lenstra's algorithm \cite[p.~538]{lenstra} therefore decides feasibility in time polynomial in \(\log(2D)\), and binary search over \(1\le m\le D\) computes \(M(L)\) within the same time bound. \qed
\end{proof}

For \(L=D\Z\times\Z\), Corollary~\ref{co:index} gives \(M(L)=D\), so listing the vertices of an optimal graph already requires \(D\) entries.

For a partition of a finite alphabet into disjoint nonempty blocks, the coordinate direct sum \(\bigoplus_{i=1}^r L_i\) consists of all integer tuples whose \(i\)-th block belongs to \(L_i\).

\begin{theorem}\label{th:blocks}
Let \(k_i\) be the size of block \(i\) in a partition of a finite alphabet into \(r\) nonempty blocks. Let \(L_i\le\Z^{k_i}\) be a full-rank sublattice for each \(i\). In the corresponding letter coordinates,
\[
 M\!\left(\bigoplus_{i=1}^r L_i\right)-1
 =\sum_{i=1}^r\bigl(M(L_i)-1\bigr).
\]
The same identity holds for unrestricted digraphs.
\end{theorem}
\begin{proof}
The phase projection in Theorem~\ref{th:main} extends verbatim to arbitrary finite alphabets: each letter translates by its corresponding coordinate unit vector, existing closed walks survive projection, and new closed walks have zero phase. Consequently \(\mn=M\) in every dimension.

Each factor admits a finite realization: its coordinate Cayley digraph realizes it by the argument of Proposition~\ref{pr:complete}, using the index times the all-ones vector in place of \((D,D)\). For the upper bound, take an optimal partial deterministic realization of each \(L_i\) and identify one vertex from each into a single root vertex. Because the alphabet blocks are disjoint, determinism is preserved. Cutting any closed walk through the root at successive visits decomposes it into single-block returns, while any closed walk avoiding the root remains within a single block. The cycle lattice of the resulting wedge is therefore \(\bigoplus_{i=1}^r L_i\), which gives the upper bound.

For the lower bound, let \(H\) be the phase image of an arbitrary realization of \(\bigoplus_{i=1}^r L_i\), with vertex set \(U\subseteq\prod_{i=1}^r Q_i\), where \(Q_i=\Z^{k_i}/L_i\). For each \(i\), project \(H\) onto the \(i\)-th factor \(U_i\subseteq Q_i\) by retaining edges labelled from block \(i\) and contracting all other steps (which preserve the \(i\)-th phase coordinate). Path projection ensures that the resulting digraph on \(U_i\) is strongly connected. Its cycle lattice contains \(L_i\), since the projections of all closed-walk counts of \(H\) generate \(L_i\), and is contained in \(L_i\), since every projected closed walk has zero displacement in \(Q_i\). Hence \(|U_i|\ge M(L_i)\).

Now choose a spanning tree of the underlying undirected graph of \(H\) and color each tree edge by its alphabet block. Because edges of colors other than \(i\) project to stationary steps on \(U_i\), the tree edges of color \(i\) project to a connected spanning subgraph of \(U_i\), so color \(i\) must appear at least \(|U_i|-1\) times. Summing over \(i=1,\ldots,r\) gives \(|U|-1\ge\sum_{i=1}^r(|U_i|-1)\), completing the proof. \qed
\end{proof}

In particular, \(M(\bigoplus_{i=1}^r d_i\Z)=1+\sum_{i=1}^r(d_i-1)\), since a unary period \(d_i\) requires a cycle of length at least \(d_i\) and a single \(d_i\)-cycle suffices. Blocks of size two contribute their rectangle values from Theorem~\ref{th:main}.

\section{Disjoint cycles and shared paths}\label{sec:packing}

\subsection{The cost of keeping cycles disjoint}

In the optimal theta construction of Theorem~\ref{th:main}, the two generating cycles necessarily share a vertex. Indeed, as we now show, forcing cycles to be pairwise vertex-disjoint incurs a precise vertex penalty governed by their lengths. Throughout this section, every cycle is simple, so its length equals its number of vertices; a loop has length one.

\begin{theorem}\label{th:packing}
Let \(G\) be a strongly connected binary labelled digraph with cycle lattice \(L\). Suppose that \(G\) contains \(k\ge2\) pairwise vertex-disjoint cycles of length at least \(g\ge1\), two of whose count vectors are linearly independent. Then
\[
 |V(G)|-M(L)\ \ge\ 1+(k-2)g.
\]
Moreover, for every pair of integers \(k\ge2\) and \(g\ge1\), equality is attained by a partial deterministic digraph whose shortest directed cycle has length \(g\).
\end{theorem}

\begin{proof}
Relabel the cycles so that the first two have independent counts \(u,v\). Write \(\ell_i\) for the length of cycle \(i\), and set \(s=\min(u_x,v_x)+\min(u_y,v_y)\). The axis-aligned rectangle spanned by \(0,u,v\) has semiperimeter
\[
 \max(u_x,v_x)+\max(u_y,v_y)=\ell_1+\ell_2-s.
\]
Since this rectangle contains two independent points of \(L\), Theorem~\ref{th:main} implies \(M(L)\le \ell_1+\ell_2-s-1\). On the other hand, vertex-disjointness yields \(|V(G)|\ge\sum_{i=1}^k\ell_i\). Subtracting these two inequalities gives the refined bound
\[
 |V(G)|-M(L)\ \ge\ 1+s+\sum_{i=3}^k\ell_i.
\]
As \(s\ge0\) and \(\ell_i\ge g\) for \(3\le i\le k\), the claimed lower bound follows.

To attain the bound, fix an integer \(m\ge\max\{g,k-1\}\). Start with a monochromatic \(b\)-cycle \(B_j\xrightarrow{b}B_{j+1}\) (with \(j\) taken modulo \(m\)). For each \(j=0,\ldots,k-2\), introduce a vertex-disjoint monochromatic \(a\)-cycle \(A_{j,i}\xrightarrow{a}A_{j,i+1}\) (with \(i\) taken modulo \(g\)). Because \(m\ge k-1\), the attachment vertices \(B_j\) are distinct, and we add the \(2(k-1)\) connecting edges
\[
 B_j\xrightarrow{a} A_{j,1},\qquad A_{j,0}\xrightarrow{b} B_{j+1}\qquad (0\le j\le k-2).
\]
Every \(a\)-cycle can be entered from and exited back to the \(b\)-cycle, so the digraph is strongly connected; furthermore, each added outgoing edge carries the label opposite to the existing outgoing edge at its tail, so the digraph is partial deterministic. Assigning phase \((0,j)\) to \(B_j\) and \((i,j)\) to \(A_{j,i}\) in \((\Z/g\Z)\times(\Z/m\Z)\), every edge---including all connecting edges---advances the phase by its corresponding coordinate unit vector. Hence every closed-walk count lies in \(g\Z\times m\Z\). The monochromatic cycles supply \((g,0)\) and \((0,m)\), so \(L=g\Z\times m\Z\). Each nonempty closed walk has a positive coordinate divisible by \(g\) in the first position or by \(m\) in the second, so its length is at least \(\min(g,m)=g\). Each \(a\)-cycle attains this length. By Corollary~\ref{co:index}, \(M(L)=g+m-1\). The \(k\) monochromatic cycles are pairwise vertex-disjoint and exhaust all \(m+(k-1)g\) vertices, giving \(|V(G)|-M(L)=1+(k-2)g\). \qed
\end{proof}

\begin{figure}[tbp]
\centering
\begin{tikzpicture}[cycle arrows,
 pt/.style={circle,draw,fill=white,inner sep=1pt,minimum size=5.2mm,font=\scriptsize},
 cyc/.style={->,semithick},link/.style={->,densely dashed,semithick},
 tag/.style={font=\scriptsize,inner sep=1pt,fill=white}]
\foreach \j/\x in {0/0,1/2.9,2/5.8,3/8.1}{\node[pt] (b\j) at (\x,0) {$B_{\j}$};}
\foreach \u/\v in {0/1,1/2,2/3}{\draw[cyc] (b\u)--(b\v);}
\draw[cyc] (b3.south) .. controls +(0,-.72) and +(0,-.72) .. (b0.south);
\foreach \j/\xx/\next in {0/0/1,1/2.9/2}{
 \begin{scope}[xshift=\xx cm]
 \node[pt] (s\j) at (.48,.88) {$1$}; \node[pt] (f\j) at (2.42,.88) {$0$};
 \node[pt] (t\j) at (1.45,1.65) {$2$};
 \draw[cyc] (f\j)--(s\j); \draw[cyc] (s\j)--(t\j); \draw[cyc] (t\j)--(f\j);
 \draw[link] (b\j)--node[tag,left]{$a$}(s\j);
 \draw[link] (f\j)--node[tag,right]{$b$}(b\next);
 \end{scope}}
\end{tikzpicture}
\caption{Equality in Theorem~\ref{th:packing} with three disjoint cycles, for \(g=3\) and \(m=4\). The bottom cycle reads \(b^4\), while each upper triangle reads \(a^3\) and is labelled by its local index \(i\). The dashed arrows are the only edges between cycles. The graph has ten vertices, four more than \(M(3\Z\times4\Z)=6\).}\label{fig:family}
\end{figure}
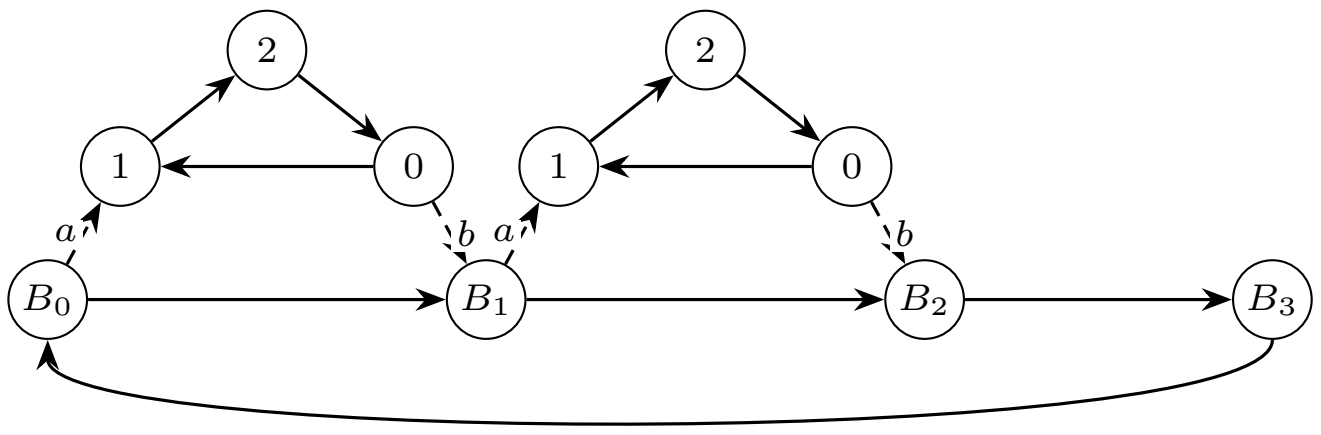

Figure~\ref{fig:split} illustrates the construction for \(k=2\) and \(g=m=3\), where identifying a single pair of vertices recovers a minimum realization; Fig.~\ref{fig:family} depicts the case \(k=3\), \(g=3\), \(m=4\). In particular, for loopless digraphs we may take \(g=2\). The excess is at least \(2k-3\), attained whenever \(m\ge\max\{2,k-1\}\).

\begin{figure}[tbp]
\centering
\begin{tikzpicture}[x=1cm,y=1cm,cycle arrows,
  pt/.style={circle,draw,fill=white,inner sep=1.2pt,minimum size=5.4mm,font=\scriptsize},
  cyc/.style={->,semithick},link/.style={->,densely dashed,semithick},
  tag/.style={font=\scriptsize,inner sep=1pt,fill=white}]
\node[pt] (a0) at (1,.78) {$A_0$}; \node[pt] (a1) at (1,-.78) {$A_1$};
\node[pt] (a2) at (-.13,0) {$A_2$};
\node[pt] (b0) at (2.86,-.78) {$B_0$}; \node[pt] (b1) at (2.86,.78) {$B_1$};
\node[pt] (b2) at (3.99,0) {$B_2$};
\foreach \u/\v in {a0/a1,a1/a2,a2/a0,b0/b1,b1/b2,b2/b0}
  {\draw[cyc] (\u)--(\v);}
\draw[link] (a0)--node[tag,above=1.5pt]{$b$}(b1);
\draw[link] (b0)--node[tag,below=1.5pt]{$a$}(a1);
\node[font=\small] at (4.85,0) {$\longrightarrow$};
\node[pt,fill=black!9] (o) at (7.54,0) {$O$};
\node[pt] (i) at (6.16,-.78) {$A_1$}; \node[pt] (j) at (6.16,.78) {$A_2$};
\node[pt] (u) at (8.92,.78) {$B_1$}; \node[pt] (v) at (8.92,-.78) {$B_2$};
\foreach \aa/\bb in {o/i,i/j,j/o,o/u,u/v,v/o}{\draw[cyc](\aa)--(\bb);}
\end{tikzpicture}
\caption{Six vertices and their five phases for \(L=3\Z^2\). Each \(A\)-triangle reads \(aaa\) and each \(B\)-triangle reads \(bbb\). On the left, the dashed arrows connect the disjoint triangles. On the right, the phase map identifies \(A_0\) and \(B_0\) as \(O\), causing the dashed arrows to coincide with existing edges so that the two cycles meet at \(O\).}\label{fig:split}
\end{figure}
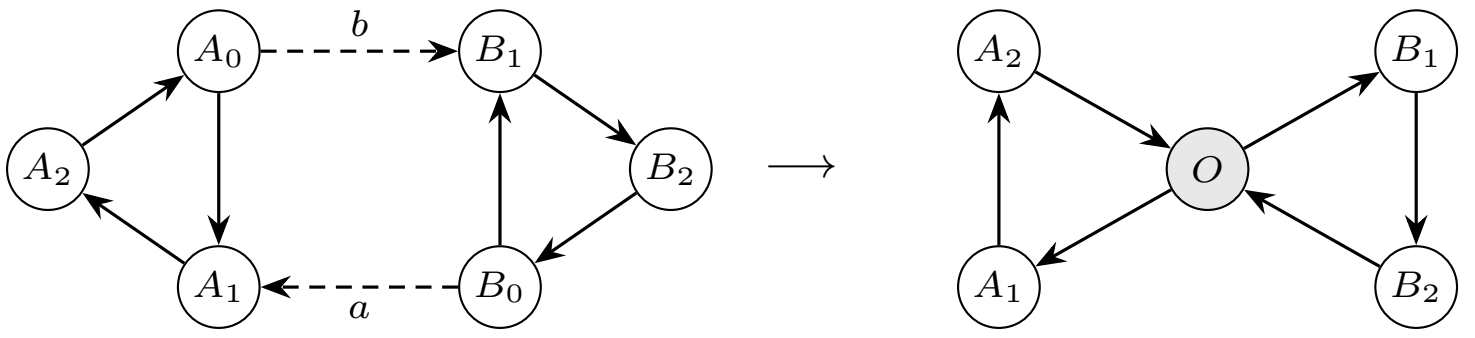

For \(g=1\), an alternative proof directly illuminates the vertex collapse under the phase map \(\phi\) of Theorem~\ref{th:main}. Let \(T_i\) be the vertex set of cycle \(i\), ordered so that the first two cycles have independent counts, and put \(U_j=\bigcup_{i=1}^j T_i\). For \(j=2\), the grid-crossing argument in Theorem~\ref{th:main} gives \(\phi(T_j)\cap\phi(U_{j-1})\ne\varnothing\); for \(j>2\), the same holds because the nonzero count of cycle \(j\) is independent of at least one of the first two counts (the grid-crossing argument relies only on periodicity and thus applies even if a projected walk repeats vertices). Setting \(d_j=|U_j|-|\phi(U_j)|\), vertex-disjointness in \(G\) implies
\[
 d_j-d_{j-1}=|T_j|-|\phi(T_j)\setminus\phi(U_{j-1})|\ge1 \qquad (2\le j\le k).
\]
Since \(d_1\ge0\), induction yields \(d_k\ge k-1\). Each vertex outside \(U_k\) contributes at most one new phase, so \(|V(G)|-|V(H)|\ge k-1\), where \(H\) is the phase quotient. Since \(H\) is a partial deterministic realization of \(L\), we have \(|V(H)|\ge M(L)\), recovering the bound for \(g=1\).

\subsection{Overlap of two basis cycles}
We next quantify how much two intersecting basis cycles can share. Let \(C_1,C_2\) be two simple cycles in a partial deterministic realization \(G\) whose counts \(u,v\) form an oppositely ordered basis of \(L(G)\), with \(u_x>v_x\) and \(v_y>u_y\). Set
\[
 d=u_x-v_x,\qquad e=v_y-u_y,\qquad D=\det(u,v).
\]
Consider the intersection of \(C_1\) and \(C_2\) as edge-labelled subgraphs, retaining isolated common vertices. Because every vertex of this intersection has indegree and outdegree at most \(1\), and a common directed cycle would force \(C_1=C_2\), the intersection is a vertex-disjoint union of directed paths (including paths of length zero). Let \(c\) denote the number of connected components of this intersection, and let \(r=(r_a,r_b)\) be the total count vector of their common edges, counted once, with coordinates indexed by their labels. Define the \emph{overlap slack}
\[
 \sigma=\min(u_x,v_x)+\min(u_y,v_y)-r_a-r_b,
\]
which measures the shortfall between the maximum number of letters potentially shareable and the number of edges actually shared. Writing \(S=u_x+v_y\), vertex inclusion--exclusion yields
\begin{equation}\label{eq:union}
 |V(C_1)\cup V(C_2)|=S+\sigma-c,
\end{equation}
since the \(c\) common paths collectively contain \(r_a+r_b+c\) vertices.

To compare positions along the two cycles, introduce the linear functional
\[
 \lambda(x,y)=\frac{ex+dy}{D},
\]
which assigns unit total advance to both basis cycles.

\begin{lemma}\label{le:order}
All vertices common to \(C_1\) and \(C_2\) appear in the same cyclic order along \(C_1\) and \(C_2\).
\end{lemma}
\begin{proof}
Since \(\lambda(u)=\lambda(v)=1\) and \(u,v\) generate \(L(G)\), we have \(\lambda(L(G))=\Z\). Fix a root vertex and assign to each vertex the value modulo \(1\) of \(\lambda\) evaluated on the count of any directed path from the root; as in Theorem~\ref{th:main}, this value is independent of the chosen path. Because \(d,e>0\), every edge contributes a strictly positive increment under \(\lambda\), and each basis cycle advances by exactly \(1\) over a full traversal. Starting from any common vertex, the partial advances along either cycle prior to return are strictly increasing in \([0,1)\). As every common vertex has the same advance modulo \(1\) along both cycles, the cyclic ordering of all common vertices is identical on \(C_1\) and \(C_2\). \qed
\end{proof}

\begin{theorem}\label{th:slack}
For every integer \(c\ge2\), the minimum possible overlap slack \(\sigma\) of two oppositely ordered basis cycles in a partial deterministic realization whose intersection has \(c\) components is \(2\lfloor c/2\rfloor\).
\end{theorem}
\begin{proof}
By Lemma~\ref{le:order}, we may index the \(c\) common paths in the same cyclic order along both cycles. For each \(i\), let \(a_i\) and \(b_i\) be the counts of the intervening arcs of \(C_1\) and \(C_2\), respectively, from the end of path \(i\) to the start of path \(i+1\) (indices modulo \(c\)). Each intervening arc is a nonempty proper subwalk of its cycle, so its \(\lambda\)-advance lies in \((0,1)\); moreover, the advances of \(a_i\) and \(b_i\) agree modulo \(1\) and hence are equal. Appending a common return path shows that \(a_i-b_i\in L(G)\). Since \(\lambda(mu+nv)=m+n\) for integers \(m,n\), the lattice vectors in the kernel of \(\lambda\) are precisely the multiples of \(u-v\). Thus there exist integers \(k_i\) such that
\[
 a_i-b_i=k_i(d,-e),\qquad \sum_{i=1}^c k_i=1,
\]
where the sum follows by subtracting the common edge count \(r\) from both \(u\) and \(v\).

Write \(a_i=(a_{ix},a_{iy})\) and \(b_i=(b_{ix},b_{iy})\), and define
\[
 T=\sum_{i=1}^c\max(0,-k_i),\qquad z=|\{i:k_i=0\}|.
\]
Put \(h_i=\min(a_{ix},b_{ix})+\min(a_{iy},b_{iy})\). Applying the identity \(2\min(x,y)=x+y-|x-y|\) to each coordinate yields
\begin{equation}\label{eq:slack}
 \sigma=\sum_{i=1}^c h_i+\frac{d+e}{2}\biggl(\sum_{i=1}^c|k_i|-1\biggr)
 =\sum_{i=1}^c h_i+(d+e)T.
\end{equation}
Whenever \(k_i=0\), we have \(a_i=b_i\); however, the initial edges of the two arcs must carry distinct labels, since partial determinism would otherwise extend the \(i\)-th common path by one edge. Hence both coordinates of \(a_i=b_i\) are strictly positive, which forces \(h_i\ge2\). Using \(d,e\ge1\) in \eqref{eq:slack}, we obtain
\[
 \sigma\ge 2z+2T,\qquad 1+2T=\sum_{i=1}^c|k_i|\ge c-z.
\]
If \(c\) is odd, these inequalities immediately give \(\sigma\ge c-1\). If \(c\) is even and \(z>0\), they give \(\sigma\ge c+z-1\ge c\). If \(c\) is even and \(z=0\), then the odd integer \(1+2T\) must be at least \(c+1\), again yielding \(\sigma\ge c\). Thus \(\sigma\ge 2\lfloor c/2\rfloor\) in all cases.

To prove sharpness, first let \(c=2k+1\). On a directed \(c\)-cycle, place two parallel edges---one labelled \(a\) and one labelled \(b\)---between each consecutive pair of vertices. Let \(C_1\) select \(a\) on \(k+1\) steps and \(b\) on the remaining \(k\) steps, and let \(C_2\) make the opposite label choice at every step. Their counts \((k+1,k)\) and \((k,k+1)\) have determinant \(2k+1\). Since every closed walk has length divisible by \(2k+1\), the cycle lattice lies in \(\{(x,y)\in\Z^2:x+y\equiv0\pmod{2k+1}\}\), which also has index \(2k+1\). The two counts therefore form an oppositely ordered basis, the \(c\) common vertices are isolated, and \(\sigma=2k\).

Now let \(c=2k\). Arrange \(2k\) common vertices in cyclic order. In one gap, insert two internally disjoint paths of length \(2\) reading \(ab\) and \(ba\), using two distinct private internal vertices; in each of the remaining \(2k-1\) gaps, place parallel edges labelled \(a\) and \(b\) as above. Let \(C_1\) choose \(a\) in \(k\) of the \(2k-1\) single-step gaps and \(b\) in the remaining \(k-1\), while \(C_2\) makes the opposite choices in those gaps, and let \(C_1,C_2\) traverse the two distinct paths in the exceptional gap. Each cycle has length \(2k+1\), and their counts are again \((k+1,k)\) and \((k,k+1)\). The same index argument shows that they form a basis; the intersection consists of \(2k\) isolated vertices and \(\sigma=2k\). \qed
\end{proof}

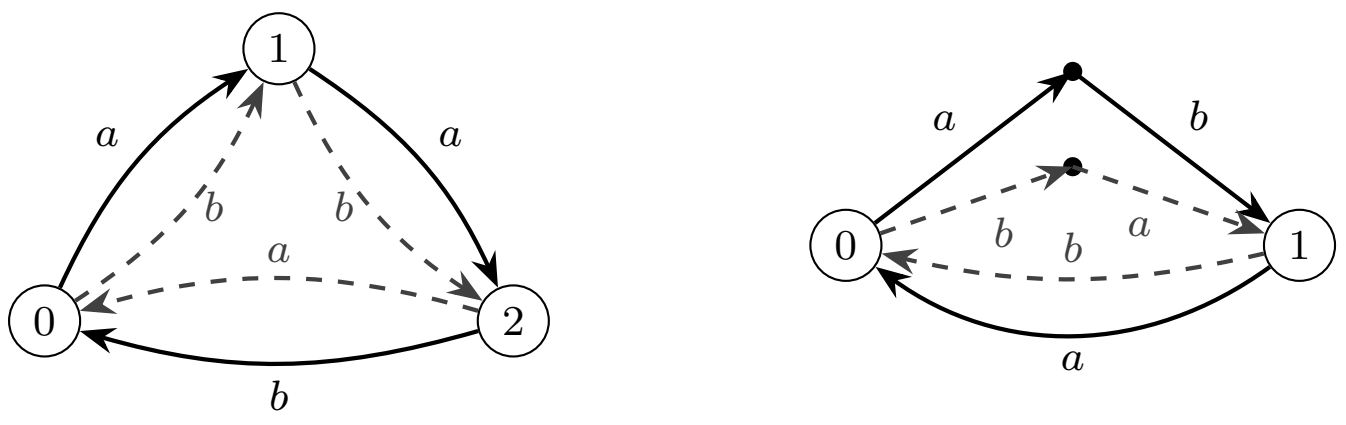
\begin{figure}[tb]
\centering
\begin{tikzpicture}[cycle arrows,font=\small,cf/.style={thick,->},cs/.style={thick,dashed,black!75,->},dot/.style={circle,draw,fill=white,minimum size=4.7mm,inner sep=0pt}]
\begin{scope}
\node[dot] (v0) at (0,0) {$0$}; \node[dot] (v1) at (1.55,1.8) {$1$}; \node[dot] (v2) at (3.1,0) {$2$};
\draw[cf] (v0) to[bend left=16] node[midway,above left] {$a$} (v1);
\draw[cs] (v0) to[bend right=16] node[midway,right] {$b$} (v1);
\draw[cf] (v1) to[bend left=16] node[midway,above right] {$a$} (v2);
\draw[cs] (v1) to[bend right=16] node[midway,left] {$b$} (v2);
\draw[cf] (v2) to[bend left=16] node[midway,below] {$b$} (v0);
\draw[cs] (v2) to[bend right=16] node[midway,above] {$a$} (v0);
\end{scope}
\begin{scope}[xshift=5.3cm,yshift=.5cm]
\node[dot] (l) at (0,0) {$0$}; \node[dot] (r) at (3.0,0) {$1$};
\coordinate (u) at (1.5,1.15); \coordinate (v) at (1.5,.52);
\fill (u) circle (1.8pt); \fill (v) circle (1.8pt);
\draw[cf] (l)--node[above left]{$a$}(u); \draw[cf] (u)--node[above right]{$b$}(r);
\draw[cs] (l)--node[below right]{$b$}(v); \draw[cs] (v)--node[below left]{$a$}(r);
\draw[cf] (r) to[bend left=36] node[midway,below]{$a$}(l);
\draw[cs] (r) to[bend left=13] node[midway,above]{$b$}(l);
\end{scope}
\end{tikzpicture}
\caption{Sharpness in Theorem~\ref{th:slack} for \(c=3\) (left) and \(c=2\) (right). Solid edges form the first cycle and dashed edges the second. Numbered vertices are common to both cycles; the two black vertices on the right are private to one cycle each. Both pairs have counts \((2,1)\) and \((1,2)\), share no edges, and attain slack two.}\label{fig:slack}
\end{figure}

Figure~\ref{fig:slack} compares the smallest odd and even extremal cases. For \(c=3\), the cycles in the odd construction may read \(aab\) and \(bba\), sharing no edges while meeting at all three vertices. Here \(u=(2,1)\), \(v=(1,2)\), and \(P=(0,0)\), giving \(\sigma=2\); equation~\eqref{eq:union} yields \(4+2-3=3\) vertices, matching Theorem~\ref{th:main} for the lattice \(x+y\equiv0\pmod3\). Since this optimum has six edges on three vertices whereas the theta realization has four, minimum-vertex realizations need not be unique even up to edge count.

\section{Random congruences of prime index}\label{sec:random}
For a prime \(p\), let \(\mathbb F_p=\Z/p\Z\) denote the finite field of order \(p\) and \(\mathbb F_p^\times\) its multiplicative group of \(p-1\) nonzero elements. Draw \(s\) uniformly at random from \(\mathbb F_p^\times\) and define
\[
 L_s=\{(x,y)\in\Z^2:x\equiv sy\pmod p\},\qquad M_s=M(L_s).
\]
Unless a different measure is indicated, all probabilities \(\Pr\), expectations \(\E\), and variances \(\Var\) below are taken with respect to this uniform choice of \(s\), and \(\log\) denotes the natural logarithm. Set \(B_s=B(L_s)=M_s+1\) (by Theorem~\ref{th:main}) and write \(Z_p=B_s/\sqrt p\) for the normalized semiperimeter. Note that the complete phase quotient has \(p\) vertices.

\begin{proposition}\label{pr:extreme}
The equality \(M_s=p\) holds if and only if \(s\equiv 1\) or \(-1\pmod p\).
\end{proposition}
\begin{proof}
Proposition~\ref{pr:complete} gives \(M_s\le p\). For \(s\equiv 1\pmod p\), any independent pair in \(L_s\cap\N^2\) must include a point off the diagonal; any such point has a coordinate at least \(p\), and a rectangle supporting two independent points has both side lengths strictly positive. For \(s\equiv -1\pmod p\), every nonzero point in \(L_s\cap\N^2\) has coordinate sum at least \(p\), whereas in a rectangle of semiperimeter at most \(p\) only the far corner attains that sum. In both cases \(B(L_s)\ge p+1\), so \(M_s=p\).

Conversely, for an integer representative \(2\le s\le p-2\), put \(j=\lceil p/s\rceil\). The vectors \((s,1)\) and \((sj-p,j)\) lie in \(L_s\cap\N^2\), have determinant \(p\), and satisfy \(0\le sj-p<s\) and \(1<j\); hence they form an oppositely ordered basis of \(L_s\). Since \(s(p-s)\ge 2(p-2)\ge p\), we have \(j\le p-s\), and therefore \(M_s\le s+j-1<p\). \qed
\end{proof}

\begin{lemma}\label{le:short}
Let \(m_s=\min_{0\ne u\in L_s}\|u\|_\infty\) be the length of a shortest nonzero vector of \(L_s\) in the maximum norm. Then \(m_s\le\sqrt p\) and
\[
 \frac{p}{m_s}\le B_s\le \frac{p}{m_s}+8m_s.
\]
\end{lemma}
\begin{proof}
By the pigeonhole principle, two distinct points in \(\{0,\ldots,\lfloor\sqrt p\rfloor\}^2\) share the same phase in \(\Z^2/L_s\); their difference lies in \(L_s\) and establishes \(m_s\le\sqrt p\). Put \(m=m_s\), and let \(u\in L_s\) attain \(\|u\|_\infty=m\). Because any common divisor of the coordinates of \(u\) is at most \(m<p\) and can be cancelled modulo \(p\), the vector \(u\) is primitive in \(\Z^2\); moreover, since \(s\not\equiv 0\pmod p\), neither coordinate of \(u\) vanishes. Any optimal rectangle for \(L_s\) contains a point \(z\in L_s\cap\N^2\) independent of \(u\), giving \(p\le|\det(u,z)|\le m(z_x+z_y)\le mB_s\), which proves the lower bound.

For the upper bound, swap coordinates if necessary (which preserves \(B_s\) and replaces \(s\) by \(s^{-1}\)) and adjust the sign of \(u\) so that \(u=(m,b)\) with \(|b|\le m\). Set \(q=p\) if \(p\ge 2m^2\), and \(q=3p\) otherwise, so that \(2m^2\le q\le p+4m^2\). By B\'ezout's identity, the line \(my-bx=q\) contains an integer point \(z'=(x,y)\). Translating \(z'\) by an integer multiple of \(u\), we may assume \(0\le x<m\), which forces \(|y-q/m|\le m\). Modulo \(p\), \(z'\) lies on the line through the nonzero vector \(u\), whence \(z'\in L_s\). The vectors \(z'\) and \(z'+u\) have determinant \(-q\ne0\) and both lie in \([0,2m]\times[0,q/m+2m]\), since \(y+b\ge q/m-2m\ge0\). Rounding the upper boundaries down to integers preserves both points, yielding a rectangle of semiperimeter at most \(q/m+4m\le p/m+8m\). \qed
\end{proof}

All asymptotics in \(p\) below are taken along primes.

\begin{theorem}\label{th:random}
For a uniformly random nonzero slope \(s\in\mathbb F_p^\times\), as \(p\to\infty\) through primes,
\[
 \E M_s=\Theta(\sqrt p),\qquad
 \Var(M_s)=\frac{12}{\pi^2}\,p\log p+O(p).
\]
\end{theorem}
\begin{proof}
Set \(R=\lfloor\sqrt p/2\rfloor\) and assume \(p\ge5\). Non-axis primitive points in \([-R,R]^2\), identified up to sign, correspond bijectively to the slopes \(s\in\mathbb F_p^\times\) with \(m_s\le R\): each such point \((a,b)\) determines \(s\equiv ab^{-1}\pmod p\), and any two candidate vectors with the same slope have determinant divisible by \(p\) and bounded in absolute value by \(2R^2<p\), so they are collinear and hence equal up to sign by primitivity. Because the shortest vector for that slope also lies in \([-R,R]^2\), the bijection preserves the maximum norm. Denoting by \(\varphi(j)\) Euler's totient (the number of integers in \(\{1,\ldots,j\}\) coprime to \(j\)), we obtain
\begin{equation}\label{eq:shell}
 |\{s\in\mathbb F_p^\times:m_s=j\}|=
 \begin{cases}
  2, & j=1,\\
  4\varphi(j), & 2\le j\le R.
 \end{cases}
\end{equation}
Indeed, for \(j\ge2\) one chooses which coordinate has magnitude \(j\), the magnitude of the other coordinate coprime to \(j\), and the relative sign, whereas for \(j=1\) only the two diagonals arise.

Write \(\mu\) for the M\"obius function. Substituting the divisor sum
\(\varphi(j)=j\sum_{d\mid j}\mu(d)/d\) from \cite[Theorem~2.3, p.~26]{apostol}
and grouping terms with \(j=dh\) gives
\[
 \sum_{j\le R}\frac{\varphi(j)}{j^2}
 =\sum_{d\le R}\frac{\mu(d)}{d^2}H_{\lfloor R/d\rfloor}
 =\frac{\log R}{\zeta(2)}+O(1),
\]
where \(H_n=\sum_{h=1}^n 1/h\). For the asymptotic equality, use
\(H_{\lfloor R/d\rfloor}=\log(R/d)+O(1)\) for \(1\le d\le R\).
The terms involving \(\log d\) have bounded total absolute value, and
\(\log R\sum_{d>R}d^{-2}=O(1)\).
The coefficient of \(\log R\) is therefore
\[
 \sum_{d\ge1}\frac{\mu(d)}{d^2}
 =\prod_{\ell\ {\rm prime}}(1-\ell^{-2})=\frac{1}{\zeta(2)}.
\]
These are the Euler products for \(\mu\) and the constant function~\(1\)
at exponent~2 \cite[Theorem~11.7 and examples, p.~231]{apostol}.
Here \(\zeta(2)=\sum_{n\ge1}n^{-2}=\pi^2/6\) by Euler's evaluation
\cite[Theorem~12.17]{apostol}. The remaining slopes with \(m_s>R\) contribute at most \(p/(R+1)^2<4\) to the unnormalized inverse-square sum. Together with \eqref{eq:shell}, this gives
\[
 \sum_{s\in\mathbb F_p^\times} m_s^{-2}
 =4\sum_{j\le R}\frac{\varphi(j)}{j^2}-2+O(1)
 =\frac{12}{\pi^2}\log p+O(1).
\]
By Lemma~\ref{le:short}, we have the pointwise estimate \(0\le B_s^2-p^2/m_s^2\le 16p+64m_s^2\le 80p\). Averaging over all \(p-1\) slopes yields \(\E[B_s^2]=(12/\pi^2)p\log p+O(p)\). Likewise, \(\sum_{s} m_s^{-1}\le 4R+p/(R+1)=O(\sqrt p)\), so \(\E B_s=O(\sqrt p)\), while Corollary~\ref{co:index} gives \(B_s\ge 2\sqrt p\). Subtracting \((\E B_s)^2=O(p)\) and using \(\Var(M_s)=\Var(B_s)\) completes the proof. \qed
\end{proof}

The same shell count yields a matching lower bound for the upper tail of \(Z_p\). For \(18\le t\le\sqrt p/4\), put \(K=\lfloor\sqrt p/(2t)\rfloor\). At least \(K^2/4\) positive pairs in \([1,K]^2\) are coprime (since \(\sum_{d\ge2}d^{-2}<3/4\)), producing distinct slopes with \(m_s\le K\) and hence \(Z_p>t\). As \(K\ge\sqrt p/(4t)\), it follows that \(\Pr(Z_p>t)\ge C/t^2\) for an absolute constant \(C>0\).

The standardized variable converges to zero in probability. Indeed, since \(M_s\ge0\),
\[
 \frac{\E|M_s-\E M_s|}{\sqrt{\Var(M_s)}}
 \le\frac{2\E M_s}{\sqrt{\Var(M_s)}}
 =O\bigl((\log p)^{-1/2}\bigr).
\]
Rescaling by \(\sqrt p\), however, yields a nondegenerate limit law.

Recall that a lattice in \(\R^2\) is \emph{unimodular} if its fundamental parallelogram has area \(1\), and that \(\mathrm{SL}_2(K)\) denotes the group of \(2\times2\) matrices of determinant \(1\) over the ring \(K\), for \(K=\R\) or \(\Z\). Right multiplication by \(\mathrm{SL}_2(\Z)\) preserves the oriented lattice generated by the columns of a basis matrix, identifying the space of unimodular lattices with \(X=\mathrm{SL}_2(\R)/\mathrm{SL}_2(\Z)\).

\begin{theorem}\label{th:limit}
Let \(\mu_X\) be the \(\mathrm{SL}_2(\R)\)-invariant (Haar) probability measure on \(X=\mathrm{SL}_2(\R)/\mathrm{SL}_2(\Z)\). For \(\Lambda\in X\), let \(b(\Lambda)\) denote the minimum semiperimeter of a closed axis-aligned real rectangle anchored at the origin that contains two linearly independent points of \(\Lambda\cap[0,\infty)^2\). If \(\Lambda\) is distributed according to \(\mu_X\), then as \(p\to\infty\) through primes, the following limits hold, where \(\xrightarrow{\mathrm d}\) denotes convergence in distribution:
\[
 Z_p\ \xrightarrow{\ \mathrm{d}\ }\ b(\Lambda),
 \qquad \frac{\E M_s}{\sqrt p}\longrightarrow\int_X b\dd\mu_X\ \in [2,\infty).
\]
\end{theorem}
\begin{proof}
Every unimodular lattice \(\Lambda\in X\) contains two linearly independent points in \((0,\infty)^2\) (obtained by approximating sufficiently large multiples of two independent interior vectors within a bounded fundamental parallelogram), and discreteness guarantees that the infimum defining \(b(\Lambda)\) is attained. By scaling, \(Z_p=b(p^{-1/2}L_s)\). The Hecke equidistribution theorem of Goldstein and Mayer \cite{hecke} (stated in lattice form in \cite[Theorem~3]{kaminaga}), applied in dimension two and along prime indices, gives weak convergence to \(\mu_X\) of the uniform measures on the lattices \(p^{-1/2}L\) with \([\Z^2:L]=p\). These sublattices correspond to the \(p+1\) lines modulo \(p\), and omitting the two coordinate lines does not affect weak convergence.

In a local basis chart on \(X\), the function \(b\) is a countable infimum of Borel functions indexed by independent pairs of integer coefficient vectors (assigning infinite cost outside \([0,\infty)^2\)). If \(\Lambda\) has no nonzero point on either coordinate axis, then any minimizing pair remains strictly inside \((0,\infty)^2\) under small perturbations of \(\Lambda\), which yields upper semicontinuity of \(b\) at \(\Lambda\). Consequently, nearby minimizers lie in a common bounded box, and the inverse basis matrices remain bounded; along any subsequence attaining the lower limit, we may therefore fix both integer coefficient vectors, whose limits remain feasible and independent, establishing lower semicontinuity as well. The exceptional set of lattices possessing a nonzero point on an axis has \(\mu_X\)-measure zero, as each fixed nonzero integer coefficient vector imposes a proper hypersurface condition. The continuous mapping theorem \cite[Theorem~2.7]{billingsley} now gives \(Z_p\xrightarrow{\mathrm{d}} b(\Lambda)\), and the deterministic shift \(p^{-1/2}\) gives the same law in the limit for \(M_s/\sqrt p\). For fixed \(t\ge18\), the lower bound on the upper tail obtained above applies for all sufficiently large primes. The closed-set inequality of the portmanteau theorem \cite[Theorem~2.1]{billingsley} gives
\[
 \Pr_{\mu_X}\{b(\Lambda)\ge t\}
 \ge \limsup_{p\to\infty}\Pr\{Z_p\ge t\}\ge C/t^2.
\]
Thus the finite random variable \(b(\Lambda)\) has unbounded support and its law is nondegenerate.

Finally, for \(t\ge18\), Lemma~\ref{le:short} implies that \(Z_p>t\) forces \(m_s<2\sqrt p/t\). The box of that radius contains \(O(p/t^2)\) nonzero integer vectors (and none when the radius is below \(1\)), each determining at most one allowed slope. Hence \(\Pr(Z_p>t)=O(t^{-2})\) uniformly in \(p\). Truncating \(Z_p\) at height \(T\) alters \(\E Z_p\) by \(O(1/T)\), while the bounded truncations converge in expectation by the distributional limit. Letting \(T\to\infty\) yields the convergence of means, since \(Z_p\ge2\) and \(M_s/\sqrt p=Z_p-p^{-1/2}\). \qed
\end{proof}

\begin{figure}[tbp]
\centering
\begin{tikzpicture}[x=1.86cm,y=.88cm,font=\footnotesize,cycle arrows]
\foreach \y in {1,2,3,4}{\draw[black!15] (.57,\y)--(4.75,\y);}
\draw[->] (.57,0)--(4.86,0) node[right] {\(p\)};
\draw[->] (.57,0)--(.57,4.35);
\node[rotate=90] at (.10,2.15) {moment / \((p\log p)\)};
\foreach \y in {0,1,2,3,4}{\draw (.57,\y)--++(-.055,0) node[left]{\y};}
\foreach \x/\v in {1/10,2/10^2,3/10^3,4/10^4}{\draw (\x,0)--++(0,-.11) node[below]{\(\v\)};}
\draw[dashed,thin] (.57,1.2158542)--(4.75,1.2158542);
\node[above,fill=white,inner sep=1pt,font=\scriptsize] at (2.05,1.2158542){\(12/\pi^2\)};
\fill (0.6989700,3.7901431) circle[radius=1.35pt];
\draw[fill=white] (0.6989700,0.0310667) circle[radius=1.5pt];
\fill (0.8450980,3.3281036) circle[radius=1.35pt];
\draw[fill=white] (0.8450980,0.0652569) circle[radius=1.5pt];
\fill (1.0413927,2.9192267) circle[radius=1.35pt];
\draw[fill=white] (1.0413927,0.1152526) circle[radius=1.5pt];
\fill (1.1139434,2.7490921) circle[radius=1.35pt];
\draw[fill=white] (1.1139434,0.1366215) circle[radius=1.5pt];
\fill (1.2304489,2.6186230) circle[radius=1.35pt];
\draw[fill=white] (1.2304489,0.1631773) circle[radius=1.5pt];
\fill (1.2787536,2.5620703) circle[radius=1.35pt];
\draw[fill=white] (1.2787536,0.1752182) circle[radius=1.5pt];
\fill (1.3617278,2.5035295) circle[radius=1.35pt];
\draw[fill=white] (1.3617278,0.1927555) circle[radius=1.5pt];
\fill (1.4623980,2.3735992) circle[radius=1.35pt];
\draw[fill=white] (1.4623980,0.2205363) circle[radius=1.5pt];
\fill (1.4913617,2.3402997) circle[radius=1.35pt];
\draw[fill=white] (1.4913617,0.2267028) circle[radius=1.5pt];
\fill (1.5682017,2.2446107) circle[radius=1.35pt];
\draw[fill=white] (1.5682017,0.2478303) circle[radius=1.5pt];
\fill (1.6127839,2.2448966) circle[radius=1.35pt];
\draw[fill=white] (1.6127839,0.2564095) circle[radius=1.5pt];
\fill (1.6334685,2.2412205) circle[radius=1.35pt];
\draw[fill=white] (1.6334685,0.2590446) circle[radius=1.5pt];
\fill (1.6720979,2.2349780) circle[radius=1.35pt];
\draw[fill=white] (1.6720979,0.2673249) circle[radius=1.5pt];
\fill (1.7242759,2.1940895) circle[radius=1.35pt];
\draw[fill=white] (1.7242759,0.2785284) circle[radius=1.5pt];
\fill (1.7708520,2.1713808) circle[radius=1.35pt];
\draw[fill=white] (1.7708520,0.2897936) circle[radius=1.5pt];
\fill (1.7853298,2.1321567) circle[radius=1.35pt];
\draw[fill=white] (1.7853298,0.2944964) circle[radius=1.5pt];
\fill (1.8260748,2.1045384) circle[radius=1.35pt];
\draw[fill=white] (1.8260748,0.3050846) circle[radius=1.5pt];
\fill (1.8512583,2.1220155) circle[radius=1.35pt];
\draw[fill=white] (1.8512583,0.3083778) circle[radius=1.5pt];
\fill (1.8633229,2.0807378) circle[radius=1.35pt];
\draw[fill=white] (1.8633229,0.3133366) circle[radius=1.5pt];
\fill (1.8976271,2.0863497) circle[radius=1.35pt];
\draw[fill=white] (1.8976271,0.3200284) circle[radius=1.5pt];
\fill (1.9190781,2.0792276) circle[radius=1.35pt];
\draw[fill=white] (1.9190781,0.3248926) circle[radius=1.5pt];
\fill (1.9493900,2.0723652) circle[radius=1.35pt];
\draw[fill=white] (1.9493900,0.3295367) circle[radius=1.5pt];
\fill (1.9867717,2.0417980) circle[radius=1.35pt];
\draw[fill=white] (1.9867717,0.3375610) circle[radius=1.5pt];
\fill (2.0043214,2.0296607) circle[radius=1.35pt];
\draw[fill=white] (2.0043214,0.3429073) circle[radius=1.5pt];
\fill (2.7015680,1.8220665) circle[radius=1.35pt];
\draw[fill=white] (2.7015680,0.4839974) circle[radius=1.5pt];
\fill (3.0038912,1.7503487) circle[radius=1.35pt];
\draw[fill=white] (3.0038912,0.5366197) circle[radius=1.5pt];
\fill (3.4772660,1.6814477) circle[radius=1.35pt];
\draw[fill=white] (3.4772660,0.6108553) circle[radius=1.5pt];
\fill (4.0003039,1.6187901) circle[radius=1.35pt];
\draw[fill=white] (4.0003039,0.6804617) circle[radius=1.5pt];
\fill (4.3009648,1.5899241) circle[radius=1.35pt];
\draw[fill=white] (4.3009648,0.7142040) circle[radius=1.5pt];
\fill (4.6021577,1.5649003) circle[radius=1.35pt];
\draw[fill=white] (4.6021577,0.7447141) circle[radius=1.5pt];
\fill (3.03,3.72) circle[radius=1.35pt];
\node[right] at (3.10,3.72){\(\E[B_s^2]/(p\log p)\)};
\draw[fill=white] (3.03,3.12) circle[radius=1.5pt];
\node[right] at (3.10,3.12){\(\Var(M_s)/(p\log p)\)};
\end{tikzpicture}
\caption{Normalized second moment (filled circles) and variance (open circles) for the 30 enumerated primes. The horizontal scale is logarithmic. The dashed line marks the common asymptotic coefficient from Theorem~\ref{th:random}; the vertical distance between the two points at each prime is the centering term \((\E B_s)^2/(p\log p)\).}\label{fig:moments}
\end{figure}
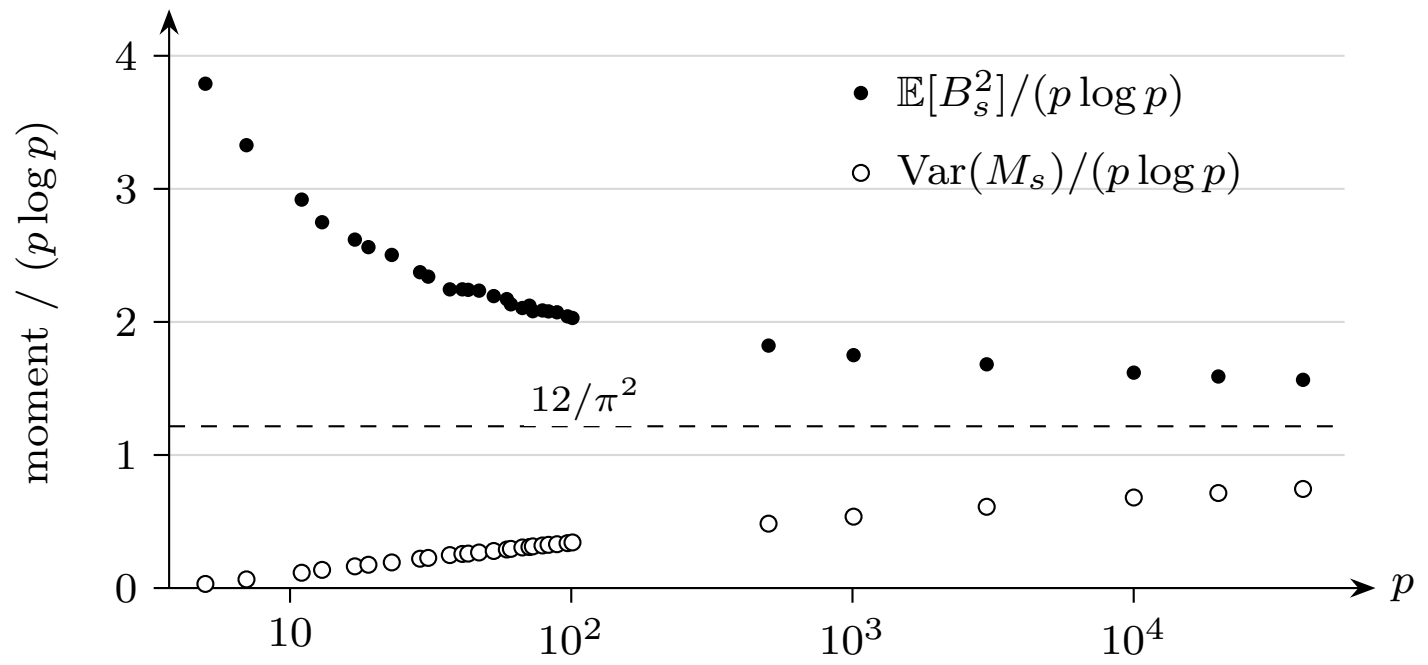

\paragraph{Numerical verification.}
An exhaustive C++ enumeration covers all \(75{,}652\) nonzero slopes across \(30\) primes from \(5\) to \(40{,}009\), verifying Lemma~\ref{le:short} and the shell counts \eqref{eq:shell}, supplemented by an independent enumeration of feasible point pairs for the listed primes up to \(101\). Table~\ref{tab:moments} separates the two terms contributing to the variance using the exact sums archived in Online Resource~1. Figure~\ref{fig:moments} compares the normalized second moment and variance across all 30 primes. At \(p=40{,}009\), centering subtracts approximately \(0.82\) from the normalized second moment; both terms remain quantitatively significant at this range, consistent with the \(O(1/\log p)\) error permitted in the variance ratio by Theorem~\ref{th:random}.

\begin{table}[htb]
\caption{Moments divided by \(p\log p\) for the six largest primes in the enumeration. The last column is the difference of the two preceding columns, computed before rounding.}\label{tab:moments}
\centering
\begin{tabular}{r r r r}
\toprule
\(p\) & \(\E[B_s^2]/(p\log p)\) & \((\E B_s)^2/(p\log p)\) & \(\Var(M_s)/(p\log p)\)\\
\midrule
503 & 1.82207 & 1.33807 & 0.48400\\
1,009 & 1.75035 & 1.21373 & 0.53662\\
3,001 & 1.68145 & 1.07059 & 0.61086\\
10,007 & 1.61879 & 0.93833 & 0.68046\\
19,997 & 1.58992 & 0.87572 & 0.71420\\
40,009 & 1.56490 & 0.82019 & 0.74471\\
\bottomrule
\end{tabular}
\end{table}

\section{Concluding remarks}\label{sec:concluding}
Several natural questions remain open. First, in the prime-index model of Section~\ref{sec:random}, it would be interesting to determine the exact value of the Haar mean \(\int_X b\dd\mu_X\) in Theorem~\ref{th:limit}, as well as whether \(\Var(M_s)-(12/\pi^2)p\log p\) admits an asymptotic linear term \(c\,p\) as \(p\to\infty\). Second, for indecomposable full-rank sublattices over alphabets of size \(k\ge3\), our phase-quotient reduction still shows \(\mn(L)=M(L)\), but the grid-crossing argument in Theorem~\ref{th:main} no longer follows from a one-dimensional discrete intermediate value argument; characterizing \(M(L)\) in higher dimensions remains an open problem.

\section*{Statements and Declarations}
\paragraph{Funding.} The research was supported by the
Hebei Province introduced overseas student support projects
 (Grant Nos.C20190365).
\paragraph{Competing interests.} No competing interests have been reported.
\paragraph{Data, code and formal proofs.} Online Resource~1 contains the exact data, the C++ enumeration, and the Lean~4.19.0 proofs of the main theorem and its optimal theta construction. The supplement lists the verification status of each result in \texttt{lean/COVERAGE.md}.
\paragraph{Use of generative AI.} The first author used OpenAI's ChatGPT to explore proof strategies, draft and revise the exposition, and assist with computations, Lean proofs and TikZ figures. The authors take responsibility for the choice of arguments, the interpretation of the results and the final manuscript, including its proofs and computational materials.


\begin{thebibliography}{20}
\bibitem{apostol} Apostol, T.M.: Introduction to Analytic Number Theory. Springer, New York (1976). \rdoi{10.1007/978-1-4757-5579-4}
\bibitem{flows} Bacher, R., de la Harpe, P., Nagnibeda, T.: The lattice of integral flows and the lattice of integral cuts on a finite graph. Bull. Soc. Math. France \textbf{125}, 167--198 (1997). \rdoi{10.24033/bsmf.2303}
\bibitem{beck} Beck, M., Robins, S.: Computing the Continuous Discretely: Integer-Point Enumeration in Polyhedra, 2nd edn. Springer, New York (2015). \rdoi{10.1007/978-1-4939-2969-6}
\bibitem{billingsley} Billingsley, P.: Convergence of Probability Measures, 2nd edn. Wiley, New York (1999). \rdoi{10.1002/9780470316962}
\bibitem{cadilhac} Cadilhac, M., Krebs, A., McKenzie, P.: The algebraic theory of Parikh automata. Theory Comput. Syst. \textbf{62}, 1241--1268 (2018). \rdoi{10.1007/s00224-017-9817-2}
\bibitem{routing} Cai, J.-Y., Havas, G., Mans, B., Nerurkar, A., Seifert, J.-P., Shparlinski, I.: On routing in circulant graphs. In: Computing and Combinatorics, Lecture Notes in Computer Science, vol. 1627, pp. 360--369 (1999). \rdoi{10.1007/3-540-48686-0_36}
\bibitem{geometry} Champarnaud, J.-M., Dubernard, J.-P., Guingne, F., Jeanne, H.: Geometrical regular languages and linear Diophantine equations: the strongly connected case. Theor. Comput. Sci. \textbf{449}, 54--63 (2012). \rdoi{10.1016/j.tcs.2012.04.012}
\bibitem{orbits} Echols, I., Harrison, J., Hudgins, T.: Periodic orbits on 2-regular circulant digraphs. Graphs Combin. \textbf{41}, 67 (2025). \rdoi{10.1007/s00373-025-02931-8}
\bibitem{farb} Farb, B., Margalit, D.: A Primer on Mapping Class Groups. Princeton University Press, Princeton (2012). \rdoi{10.1515/9781400839049}
\bibitem{hecke} Goldstein, D., Mayer, A.: On the equidistribution of Hecke points. Forum Math. \textbf{15}, 165--189 (2003). \rdoi{10.1515/form.2003.009}
\bibitem{hoffmann} Hoffmann, S.: State complexity bounds for the commutative closure of group languages. J. Autom. Lang. Comb. \textbf{28}, 27--57 (2023). \rdoi{10.25596/JALC-2023-027}
\bibitem{kaminaga} Kaminaga, M.: Upper bound for the coefficients of the shortest vector of random lattice. IEICE Trans. Fundam. Electron. Commun. Comput. Sci. \textbf{E106-A}, 1585--1588 (2023). \rdoi{10.1587/transfun.2023EAL2032}
\bibitem{klaedtke} Klaedtke, F., Rue\ss{}, H.: Monadic second-order logics with cardinalities. In: Automata, Languages and Programming, Lecture Notes in Computer Science, vol. 2719, pp. 681--696 (2003). \rdoi{10.1007/3-540-45061-0_54}
\bibitem{lavado} Lavado, G.J., Pighizzini, G., Seki, S.: Converting nondeterministic automata and context-free grammars into Parikh equivalent one-way and two-way deterministic automata. Inf. Comput. \textbf{228--229}, 1--15 (2013). \rdoi{10.1016/j.ic.2013.06.003}
\bibitem{lenstra} Lenstra, H.W., Jr.: Integer programming with a fixed number of variables. Math. Oper. Res. \textbf{8}, 538--548 (1983). \rdoi{10.1287/moor.8.4.538}
\bibitem{diameters} Marklof, J., Str\"ombergsson, A.: Diameters of random circulant graphs. Combinatorica \textbf{33}, 429--466 (2013). \rdoi{10.1007/s00493-013-2820-6}
\bibitem{spines} Martelli, B., Novaga, M., Pluda, A., Riolo, S.: Spines of minimal length. Ann. Sc. Norm. Super. Pisa Cl. Sci. (5) \textbf{17}, 1067--1090 (2017). \rdoi{10.2422/2036-2145.201511_003}
\bibitem{nerode} Nerode, A.: Linear automaton transformations. Proc. Amer. Math. Soc. \textbf{9}, 541--544 (1958). \rdoi{10.1090/S0002-9939-1958-0135681-9}
\bibitem{parikh} Parikh, R.J.: On context-free languages. J. ACM \textbf{13}, 570--581 (1966). \rdoi{10.1145/321356.321364}
\bibitem{schrijver} Schrijver, A.: Graphs on the torus and geometry of numbers. J. Combin. Theory Ser. B \textbf{58}, 147--158 (1993). \rdoi{10.1006/jctb.1993.1033}
\end{thebibliography}
\end{document}